\documentclass{amsart}
\usepackage{amssymb,psfig}

\newcommand{\bZ}{\mbox{${\mathbb Z}$}}

\newcommand{\bC}{\mbox{${\mathbb C}$}}

\newcommand{\setp}[2]{\mbox{$\{#1\::\:#2\}$}}

\newcommand{\case}[5]{#1#2 \left\{ \begin{array}{ll} #3 &\mbox{if}\spa #4 \\ #5 &\mbox{otherwise\,.} \end{array} \right.}
\newcommand{\casec}[5]{#1#2 \left\{ \begin{array}{ll} #3 &\mbox{if}\spa #4 \\ #5 &\mbox{otherwise\,,} \end{array} \right.}
\newcommand{\casecc}[6]{#1#2 \left\{ \begin{array}{ll} #3 &\mbox{if}\spa #4 \\ #5 &\mbox{if}\spa #6\,, \end{array} \right.}

\newcommand{\casex}[3]{ \left\{ \begin{array}{ll} #1 &\mbox{if}\spa #2 \\ #3 &\mbox{otherwise\,} \end{array} \right.}
\newcommand{\casexsc}[3]{ \left\{ \begin{array}{ll} #1 &\mbox{if}\spa #2 \\ #3 &\mbox{otherwise\,;} \end{array} \right.}

\newcommand{\stacksum}[2]{\sum_{\begin{array}{c}\vspace{-5mm}\;\\ \vspace{-1mm}\scriptstyle{#1}\\ \scriptstyle{#2}\end{array}} }

\newcommand{\eproof}{\hfill\mbox{$\Box$}}
\newcommand{\spa}{\;\:}

\newcommand{\fg}{{\mathfrak G}}
\newcommand{\fs}{{\mathfrak S}}

\newcommand{\s}{\sigma}

\newtheorem{conj}[equation]{Conjecture}

\newtheorem{thm}[equation]{Theorem}
\newtheorem{rem}[equation]{Remarks}
\newtheorem{cor}[equation]{Corollary}
\newtheorem{lem}[equation]{Lemma}
\newtheorem{exa}[equation]{Example}

\makeatletter
\let\choose\@@choose

\makeatother

\makeatletter\@addtoreset{equation}{section}\makeatother

\begin{document}
\bibliographystyle{plain}
\title[$K^0(Fl_n)$ and the Fomin-Kirillov Quadratic Algebra]{The $K$-theory of the Flag Variety and the Fomin-Kirillov Quadratic Algebra}
\author{Cristian Lenart}
\address{Department of Mathematics and Statistics, State University of New York, Albany, NY 12222}
\email{lenart@csc.albany.edu}

\thanks{The author was supported by SUNY Albany Faculty Research Award 1032354}

\begin{abstract}
We propose a new approach to the multiplication of Schubert classes in the $K$-theory of the flag variety. This extends the work of Fomin and Kirillov in the cohomology case, and is based on the quadratic algebra defined by them. More precisely, we define $K$-theoretic versions of the Dunkl elements considered by Fomin and Kirillov, show that they commute, and use them to describe the structure constants of the $K$-theory of the flag variety with respect to its basis of Schubert classes.
\end{abstract}

\vspace{-2cm}

\maketitle

\vspace{-1.1cm}

\section{Introduction}

An important open problem in algebraic combinatorics is to describe combinatorially the structure constants for the cohomology of the flag variety $Fl_n$ (that is, the variety of complete flags $(0=V_0\subset V_1\subset\ldots \subset V_{n-1}\subset V_n=\bC^n)$ in $\bC^n$) with respect to its basis of Schubert classes. These structure constants are known as {\em Littlewood-Richardson coefficients}; a subset of them, consisting of the structure constants for the cohomology of a Grassmannian, are described by the classical Littlewood-Richardson rule. Fomin and Kirillov \cite{fakqad} proposed a new approach to the Littlewood-Richardson problem for $H^*(Fl_n)$ based on a certain algebra with quadratic relations that they defined. In this paper, we extend Fomin and Kirillov's approach to the $K$-theory of $Fl_n$. 

It is well-known that the integral cohomology ring $H^*(Fl_n)$ and the Grothendieck ring $K^0(Fl_n)$ are both isomorphic to $\bZ[x_1,\ldots,x_n]/I_n$, where $I_n$ is the ideal generated by symmetric polynomials in $x_1,\ldots,x_n$ with constant term 0. In the cohomology case, the elements $x_i$ are identified with the Chern classes of the dual line bundles $L_i^*$, where $L_i:=V_i/V_{i-1}$ are tautological line bundles. In the $K$-theory case, we identify $x_i$ with the $K$-theory Chern class $1-1/y_i$ of the line bundle $L_i^*$, where $y_i:=1/(1-x_i)$ represents $L_i$ in the Grothendieck ring.

One can define natural bases for $H^*(Fl_n)$ and $K^0(Fl_n)$ (over $\bZ$) based on the CW-complex structure of $Fl_n$ given by the (opposite) {\em Schubert varieties}. These are varieties $X_w$, which are indexed by permutations $w$ in $S_n$, and which have complex codimension $l(w)$ (that is, the number of inversions in $w$). More precisely, if we think of $Fl_n$ as $SL_n/B$, we let $X_w:=\overline{B^-wB/B}$, where $B$ and $B^-$ are the subgroups of $SL_n$ consisting of upper and lower triangular matrices. The  {\em Schubert} and {\em Grothendieck polynomials} indexed by $w$, which are denoted by $\fs_w(x)$ and $\fg_w(x)$, are certain polynomial representatives for the cohomology and $K$-theory classes corresponding to $X_w$. These classes, which are denoted by $\sigma_w$ and $\omega_w$, form the mentioned natural bases for $H^*(Fl_n)$ and $K^0(Fl_n)$. Schubert and Grothendieck polynomials were defined by Lascoux and Sch\"utzenberger \cite{lasps,lassha}, and were studied extensively during the last two decades \cite{fulyt,lenuac,lrsppk,macnsp,manfsp}. 

Both the Schubert polynomials $\fs_w(x)$ and the Grothendieck polynomials $\fg_w(x)$, for $w$ in $S_{\infty}$, form bases of $\bZ[x_1,x_2,\ldots]$; here $S_{\infty}:=\bigcup_n S_n$ under the usual inclusion $S_n\hookrightarrow S_{n+1}$. Hence we can write
\[\fs_u(x)\,\fs_v(x)=\sum_{w\::\:l(w)=l(u)+l(v)}c_{uv}^w\fs_w(x)\,,\spa\spa\;\; \fg_u(x)\,\fg_v(x)=\sum_{w\::\:l(w)\ge l(u)+l(v)}c_{uv}^w\fg_w(x)\,.\]
The notation is consistent since the structure constants corresponding to Schubert polynomials, which are known to be nonnegative for geometric reasons, are a subset of the structure constants corresponding to Grothendieck polynomials. 

The simplest multiplication formula for Schubert polynomials is {\em Monk's formula}, which can be stated as follows:
\begin{equation}\label{monk}
x_p\,\fs_v(x)=-\stacksum{1\le i<p}{l(vt_{ip})=l(v)+1}\fs_{vt_{ip}}(x)+\stacksum{i>p}{l(vt_{pi})=l(v)+1}\fs_{vt_{pi}}(x)\,;
\end{equation}
here $t_{ij}$ is the transposition of $i,j$, and $v$ is an arbitrary element of $S_\infty$. In fact, Monk's formula expresses the product of $\fs_v(x)$ with a Schubert polynomial indexed by an adjacent transposition, which is equivalent to (\ref{monk}). 

Similar formulas for Grothendieck polynomials were derived in \cite{lenktv}. Define the set $\varPi_p(v)$ to consist of all permutations
\[w=vt_{i_1p}\ldots t_{i_rp}t_{pi_{r+1}}\ldots t_{pi_{r+s}}\,,\]
in $S_\infty$ such that $r,s\ge 0$, $r+s\ge 1$, the length increases by precisely 1 upon multiplication by each transposition, and 
\[i_r<\ldots<i_1<p<i_{r+s}<\ldots<i_{r+1}\,.\]
Given a permutation $w$ in $\varPi_p(v)$, let $\s_p(w,v):=(-1)^{s+1}$.

\begin{thm}\label{kmonk}\cite{lenktv}.
We have
\[x_p\,\fg_v(x)=\sum_{w\in \varPi_p(v)} \s_p(w,v)\fg_w(x)\,.\]
\end{thm}

Motivated by the Littlewood-Richardson problem for Schubert polynomials, Fomin and Kirillov \cite{fakqad} defined the {\em quadratic algebra} ${\mathcal E}_n$ (of type $A_{n-1}$) as the associative algebra with generators $[ij]$ for $1\le i\ne j\le n$, which satisfy the following relations:
\begin{enumerate}
\item[(i)] $[ij]+[ji]=0\,$,
\item[(ii)] $[ij]^2=0\,$,
\item[(iii)] $[ij][jk]+[jk][ki]+[ki][ij]=0\,$,
\item[(iv)] $[ij][kl]=[kl][ij]\,$,
\end{enumerate}
for all distinct $i,j,k,l\,$. There is a natural grading on ${\mathcal E}_{n}$ given by ${\rm deg}([ij])=1$. This algebra provides a solution to the classical Yang-Baxter equation.

There is a natural representation of the quadratic algebra, called the {\em Bruhat representation}. This is a representation on the group algebra $\bZ\langle S_n\rangle$ of the symmetric group, which is identified (as a vector space) with the cohomology of the flag variety $H^*(Fl_n)$, via $w\mapsto\sigma_w$. The representation is defined by
\[\casec{[ij]\,w}{=}{wt_{ij}}{l(wt_{ij})=l(w)+1}{0}\]
for $1\le i<j\le n$. 

Fomin and Kirillov defined the {\em Dunkl elements} $\theta_p$ (in the quadratic algebra) for $p=1,\ldots,n$ by
\[\theta_p:=-\sum_{1\le i<p}[ip]+\sum_{p<k\le n}[pk]=-\sum_{i\ne p}[ip]\,.\]
These elements encode the multiplicative structure of $H^*(Fl_n)$. Indeed, with the above notation, Monk's formula (\ref{monk}) can be stated as
\[\theta_p\,w=x_p\cdot\sigma_w\,;\]
this formula is to be understood under the standard vector space isomorphisms (mentioned above and at the beginning of the section) between $\bZ\langle S_n\rangle$, $H^*(Fl_n)$, and $Z[x_1,\ldots,x_n]/I_n$. 

Fomin and Kirillov showed that the Dunkl elements commute (in the quadratic algebra). Then, based on this fact, they present a new approach to the study of the multiplication of Schubert classes in cohomology, which is outlined below. 

Consider the evaluation of Schubert polynomials at Dunkl elements
\[\fs_w(\theta):=\fs_w(\theta_1,\ldots,\theta_{n-1})\,.\]
As noted in \cite{fakqad}, given $u,\,v,\,w$ in $S_n$ with $l(w)=l(u)+l(v)$, the facts stated above imply 
\begin{equation}\label{cor1}c_{uv}^w=\langle\mbox{coefficient of $w$ in $\fs_u(\theta)\,v$}\rangle\end{equation}
in the Bruhat representation. 

Let ${\mathcal E}_{n}^{+}$ be the cone of all nonnegative linear combinations of all noncommutative monomials in the generators $[ij]$, for $i<j$. Note that by applying the defining relations of ${\mathcal E}_{n}$, this cone also contains some linear combinations of monomials with negative coefficients. 

\begin{conj}\cite{fakqad}\label{conj1} (Nonnegativity conjecture). For any $w$ in $S_n$, the evaluation $\fs_w(\theta)$ lies in ${\mathcal E}_n^+$. 
\end{conj}

In Section 2, we define $K$-theoretic Dunkl elements in ${\mathcal E}_n$ based on the $K$-theory version of Monk's formula in Theorem \ref{kmonk}. We prove that these elements still commute in ${\mathcal E}_n$, but the proof is vastly more complex than in the cohomology case, due to the exponential increase in the number of terms. In Section 3, we state the $K$-theory version of the nonnegativity Conjecture \ref{conj1}, and extend to $K$-theory the above approach to the Littlewood-Richardson problem for Schubert polynomials. We also report on a joint project with A. Yong on realizing $K^0(Fl_n)$ as the commutative subalgebra of ${\mathcal E}_n$ generated by the $K$-theoretic Dunkl elements, thus extending the similar result of Fomin and Kirillov for $H^*(Fl_n)$. 

\section{$K$-theoretic Dunkl elements}

A. Yong \cite{yonpc} suggested extending the definition of Dunkl elements in \cite{fakqad}, by constructing $K$-theoretic versions of them based on Theorem \ref{kmonk}. More precisely, we define the {\em K-theoretic Dunkl element} $\kappa_p$ in ${\mathcal E}_n$ by
\[\kappa_p:=1-(1+[p-1,p])(1+[p-2,p])\ldots(1+[1p])(1+[np])(1+[n-1,p])\ldots(1+[p+1,p])\,.\]
Clearly, the degree one component of $\kappa_p$ is $\theta_p$. 

It is useful to give an expanded version of the above definition. Given $p$ with $1\le p\le n$ and a set $A\subseteq[n]\setminus \{p\}$ (where $[n]$ denotes $\{1,\ldots,n\}$), order the elements of $A=\{i_1,\ldots,i_{r+s}\}$ such that
\[i_r<\ldots<i_1<p<i_{r+s}<\ldots<i_{r+1}\,.\]
This can be thought of as a circular order (see the figure below).

\[
\begin{array}{c}
\mbox{\psfig{file=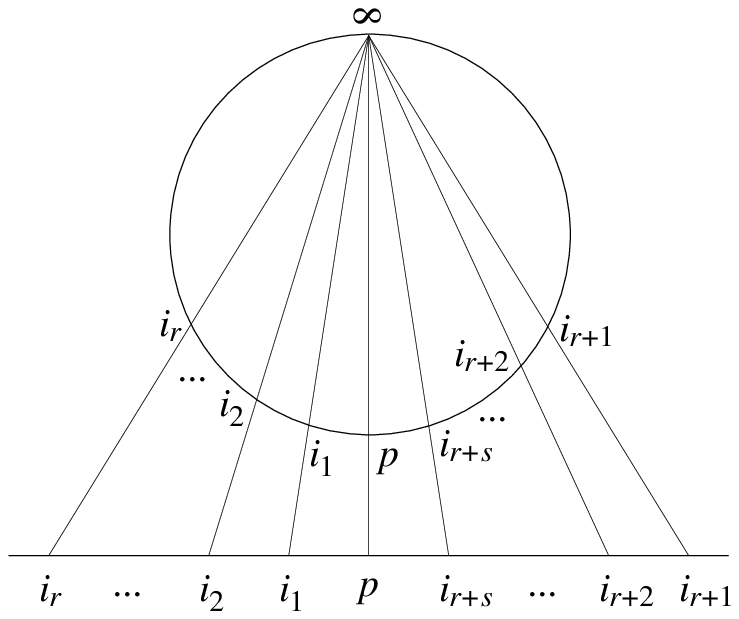}} 
\end{array}
\]

\noindent We use the predicate $C(i_1,\ldots,i_{r+s},p)$ to indicate that the elements of $A$ are in the above order. Let
\[\pi(p,A):=(-1)^{s}[i_1p]\ldots[i_rp][pi_{r+1}]\ldots[pi_{r+s}]=[i_1p]\ldots[i_{r+s}p]\]
in the quadratic algebra. Then we have
\[\kappa_p=-\sum_{\emptyset\ne A\subseteq[n]\setminus\{p\}}\pi(p,A)\,,\]
for $p=1,\ldots,n$. Note that these elements have $2^{n-1}-1$ terms with degrees between $1$ and $n-1$, compared to just $n-1$ terms of degree $1$ for the Dunkl elements in the cohomology case.

\begin{exa}{\rm 
The $K$-theoretic Dunkl elements for $n=3$ are
\[\begin{array}{lll}
\kappa_1 & = & \ \ [12]+[13]-[13][12] \\
\kappa_2 & = & -[12]+[23]+[12][23] \\
\kappa_3 & = & -[13]-[23]-[23][13]\,. 
\end{array}\]
For $n=4$, they are
\[\begin{array}{lll}
\kappa_1 & = & \ \  [12]+[13]+[14]-[13][12]-[14][12]-[14][13]+[14][13][12]
\\
\kappa_2 & = & -[12] + [23] + [24] + [12][23] + [12][24] -[24][23]
-[12][24][23]\\
\kappa_3 & = & -[13]-[23]+[34]+[13][34] -[23][13] +
[23][34]+[23][13][34]\\
\kappa_4 & = & -[14] -[24] -[34] -[24][14] -[34][14] -[34][24] -[34][24][14]\,.
\end{array}\]
}
\end{exa}

Let us identify the Bruhat representation $\bZ\langle S_n\rangle$ with $K^0(Fl_n)$ via $w\mapsto \omega_w$; also recall the ring isomorphism between $K^0(Fl_n)$ and $Z[x_1,\ldots,x_n]/I_n$ mentioned in Section 1. Given these isomorphims, we have
\[\kappa_p\,w=x_p\cdot\omega_w\,,\]
based on Theorem \ref{kmonk}. Thus, the $K$-theoretic Dunkl elements encode the multiplicative structure of the $K$-theory of the flag variety. According to the above observation, the images of Dunkl elements in the Bruhat representation commute, but this is not a faithful representation of the quadratic algebra. Nevertheless, we have the stronger result below.
\begin{thm}\label{commute} The $K$-theoretic Dunkl elements commute pairwise (in the quadratic algebra).
\end{thm}

In fact, a stronger statement holds, involving the ``restriction'' of fixed degree of the commutator to an arbitrary set 
\[\emptyset\ne X\subseteq[n]\setminus\{p,q\}\,.\]

 Given $d$ with $|X|\le d\le 2|X|+2$, we define
\begin{equation}\label{rescom}[\kappa_q,\kappa_p]_{X,d}:=\sum_{(A,B)}[\pi(p,A),\pi(q,B)]\end{equation}
where the summation ranges over all pairs of sets $(A,B)$ satisfying
\begin{itemize}
\item $X\subseteq A\cup B\subseteq X\cup\{p,q\}\,$;
\item $|A|+|B|=d\,$;
\item $A\ne\emptyset\ne B$ and $p\not\in A$, $q\not\in B\,$.
\end{itemize}
We will show that 
\[[\kappa_q,\kappa_p]_{X,d}=0\,.\] 
Furthermore, we claim that only commutations between $[ip]$ and $[iq]$ are needed (not between $[ip]$ and $[jp]$, or $[iq]$ and $[jq]$).

Before presenting the proof, we simplify the notation. Given the set $X$ above, we consider the free associative algebra of words over the alphabet $X\cup X'\cup\{p',q,*\}$, where $X':=\setp{i'}{i\in X}$. The letter $i$ corresponds to the generator $[ip]$ of the quadratic algebra, whereas $i'$ corresponds to $[iq]$. We impose the following relations on the new generators, which come from the corresponding relations (i)-(iv) defining the quadratic algebra:
\begin{enumerate}
\item $p'=-q=*\,,$
\item $i^2=i'^2=*^2=0\,,$
\item $ii'=*i-i'*$, $i'i=i*-*i'\,,$
\item $ij'=j'i$\,,
\end{enumerate}
for all $i\ne j$ in $X$. Let us denote by ${\mathcal E}_{X}$ the new algebra. Let
\[\pi(A):=i_1\ldots i_s\,,\;\;\;\pi'(B):=j_1'\ldots j_t'\,,\]
where $A=\{i_1,\ldots,i_s\}$, $B=\{j_1,\ldots,j_t\}$ are subsets of $X\cup\{p,q\}$, and $C(i_1,\ldots,i_s,p)$, $C(j_1,\ldots,j_t,q)$. Theorem \ref{commute} will follow from the Lemma below.

\begin{lem}\label{lem1} We have
\[\sum_{(A,B)}[\pi(A),\pi'(B)]=0\spa\mbox{in}\spa{\mathcal E}_{X}\,,\]
where the summation ranges over the same pairs of sets $(A,B)$ as the ones in {\rm (\ref{rescom})}. 
\end{lem}

Let us denote the above summation by $\varSigma(X,d)$. The Lemma clearly holds when $|X|=d$. On the other hand, when $d=2|X|+2$, there is just one commutator in $\varSigma(X,d)$, and its terms are both 0. Indeed, we have
\begin{equation}\label{starstar}
*i_1i_1'\ldots i_si_s'*=0\,,\;\;\;*i_1'i_1\ldots i_s'i_s*=0\,,
\end{equation}
which are easily seen by expanding all products $i_ji_j'$ and $i_j'i_j$ using the relations (3). So we assume that $|X|+1\le d\le 2|X|+1$. Let us first consider the special case when $X=\{k_1,\ldots,k_r\}$ and $C(k_1,\ldots,k_r,p,q)$. 

{\em Proof of the special case of Lemma {\rm \ref{lem1}}.} The idea is to apply the following four-step procedure:
\begin{enumerate}
\item express  $\pi(A)\pi'(B)$ and $\pi'(B)\pi(A)$ as a concatenation of two types of blocks, called {\em expansion blocks} and {\em intermediate blocks};
\item define one or two binary labels for each expansion block; 
\item expand $\varSigma(X,d)$ by replacing each expansion block with a certain sum of terms, based on the labeling;
\item group the terms obtained after expansion into classes, and show that the sum within each class is 0.
\end{enumerate}

\underline{Step 1: Preparation.} Note that the only positions in which the letter $*$ can appear in a term $\pi(A)\pi'(B)$ are the first and the last. As far as the terms $\pi'(B)\pi(A)$ are concerned, the letter $*$ can appear either at the end of the subword $\pi'(B)$ or at the beginning of the subword $\pi(A)$, but not in both positions. If $*$ appears in one of the two positions and $\pi'(B)\ne *$, respectively $\pi(A)\ne -*$, then the corresponding term $\pi'(B)\pi(A)$ can be cancelled with the one obtained by moving the $*$ from one subword to the other (according to relation (1) in the definition of ${\mathcal E}_{X}$). So we will only consider terms $\pi'(B)\pi(A)$ containing no $*$, plus the terms $*\pi(X)$ and $-\pi'(X)*$ if $d=|X|+1$. The last two terms will not be expanded, so they are excluded from the discussion in Steps 2-3 below.

We use the trivial commutation (4) in the definition of ${\mathcal E}_{X}$ in order to express $\pi(A)\pi'(B)$ as a concatenation of blocks of the form $[*]i_1i_1'\ldots i_si_s'[*]$, with $s\ge 0$, and blocks containing only letters in $(A\setminus (B\cup\{q\}))\cup(B\setminus (A\cup\{p\}))'$; here $[*]$ indicates the fact that * might or might not appear. The blocks of the first type are the expansion blocks, and the others are the intermediate ones. For $\pi'(B)\pi(A)$, the expansion blocks are of the form $j_1'j_1\ldots j_t'j_t$, with $t\ge 1$, and the intermediate blocks are defined in the same way. Note that at most one $*$ can appear in the expansion blocks of $\pi(A)\pi'(B)$, according to (\ref{starstar}). 

\underline{Step 2: Labeling.}  We consider the terms $\pi(A)\pi'(B)$ and $\pi'(B)\pi(A)$ of $\varSigma(X,d)$ separately.

\underline{Step 2.1: Labeling for $\pi(A)\pi'(B)$.} Consider such a term, and assume it has $m$ expansion blocks. If the $k$th such block is of the form $i_1i_1'\ldots i_si_s'$, we will associate to it a left label $l_k$ and a right label $r_k$, both of which are 0 or 1. For expansion blocks containing *, we define only one label; more precisely, if the $*$ is at the beginning, we define a right label $r_1$, and if the $*$ is at the end, we define a left label $l_m$. 

\[
\begin{array}{c}
\mbox{\psfig{file=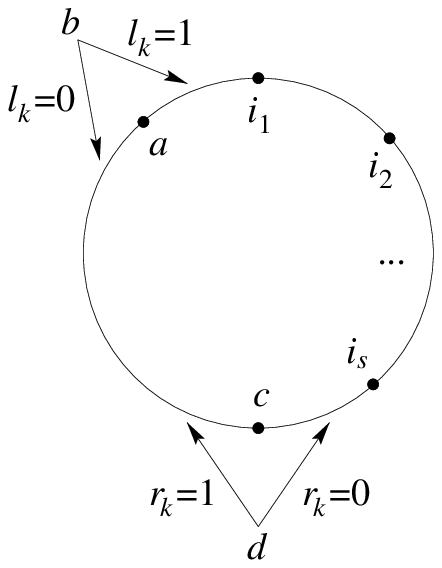}} 
\end{array}
\]

Let us now explain how $l_k$ and $r_k$ are defined (see the figure above). If the corresponding expansion block does not start with $*$, we set $l_k:=1$ in the following two cases: (1) the intermediate block preceding the $k$th expansion block ends in $ab'$ and $C(a,b,i_1)$; (2) the mentioned intermediate block consists entirely of letters in $(B\setminus (A\cup\{p\}))'$ or is empty. Otherwise, we set $l_k:=0$. Similarly, if the $k$th expansion block does not end in $*$, the right label $r_k$ is defined to be 1 in the following two cases: (1) the intermediate block after the $k$th expansion block starts with $cd'$ and $C(i_s,c,d)$; (2) the mentioned intermediate block consists entirely of letters in $A\setminus (B\cup\{q\})$ or is empty. Otherwise, we set $r_k:=0$. 

{\em Remark.} Let us explain the point of the labeling. Assume that $\pi(A)\pi'(B)=\ldots ab'i_1i_1'\ldots$. Upon expansion, we get a term of the form $\ldots ab'*i_1\ldots$. This can also be obtained by expansion from $\ldots abb'i_1\ldots$ and $\ldots b'a'ai_1\ldots$ (see below for the second expansion). But precisely one of two situations can occur: either the first word is of the form $\pi(C)\pi'(D)$ in ${\mathcal E}_X$, which happens if $l_k=1$, or the second word is of the form $\pi'(D)\pi(C)$ in ${\mathcal E}_X$, which happens if $l_k=0$; here, $(C,D)$ is a pair of sets appearing in the summation (\ref{rescom}). There is an analogous situation involving the right label.

\underline{Step 2.2: Labeling for $\pi'(B)\pi(A)$.} Similarly, we associate a left label and a right label to the expansion blocks in such a term. The left label $l_k$ of the $k$th expansion block is defined to be 1 in the following two cases: (1) the intermediate block preceding the $k$th expansion block ends in $b'a$ and $C(b,a,j_1)$; (2) the mentioned intermediate block is nonempty and consists entirely of letters in $A\setminus (B\cup\{q\})$. Otherwise, we set $l_k:=0$. Similarly, the right label $r_k$ is defined to be 1 in the following two cases: (1) the intermediate block after the $k$th expansion block starts with $d'c$ and $C(j_t,d,c)$; (2) the mentioned intermediate block is nonempty and consists entirely of letters in $(B\setminus (A\cup\{p\}))'$. Otherwise, we set $r_k:=0$. 

\underline{Step 3: Expansion.} Again, we consider the terms $\pi(A)\pi'(B)$ and $\pi'(B)\pi(A)$ of $\varSigma(X,d)$ separately.

\underline{Step 3.1: Expanding the terms $\pi(A)\pi'(B)$.} Consider such a term written in the block form described above, and assume there are $m$ expansion blocks. Let the $k$th such block be $i_1i_1'\ldots i_si_s'$. Based on the labels $l_k,\,r_k$ defined in Step 2 and the first relation (3) in the definition of ${\mathcal E}_{X}$, we express the considered block as $u_k+v_k$, where $u_k$ and $v_k$ are given by the rules below.
\begin{itemize}
\item Case 1: If $l_k=0$ and $r_k=1$ then 
\[u_k:=-i_1i_1'\ldots i_{s-1}i_{s-1}'i_s'*\,,\spa\spa v_k:=*i_1\ldots *i_s\,.\]
\item Case 2: If $l_k=1$ and $r_k=0$ then 
\[u_k:=*i_1i_2i_2'\ldots i_si_s'\,,\spa\spa v_k:=(-1)^s i_1'*\ldots i_s'*\,.\]
\item Case 3: If $l_k=1$ and $r_k=1$ then 
\[u_k:=*i_1i_2i_2'\ldots i_si_s'-i_1i_1'\ldots i_{s-1}i_{s-1}'i_s'*+*i_1i_2i_2'\ldots i_{s-1}i_{s-1}'i_s'*\,,\spa\spa v_k:=0\,;\]
if $s=1$, then the third term in $u_k$ is missing.
\item Case 4: If $l_k=0$ and $r_k=0$ then 
\[u_k:=0\,,\spa\spa v_k:=\sum_{j=0}^s (-1)^{s-j} *i_1\ldots *i_ji_{j+1}'*\ldots i_s'*\,.\]
\end{itemize}

The fact that $i_1i_1'\ldots i_si_s'=u_k+v_k$ in the four cases above can be proved by induction on $s$; the base case $s=1$ is just the first relation (3) in the definition of ${\mathcal E}_{X}$, whence all four identities are generalizations of the mentioned relation. For instance, in the last case, assuming the identity for $s-1$, we have
\begin{align*}&\sum_{j=0}^s (-1)^{s-j} *i_1\ldots *i_ji_{j+1}'*\ldots i_s'*\\
&=*i_1\ldots *i_s-\left(\sum_{j=0}^{s-1}(-1)^{s-1-j} *i_1\ldots *i_ji_{j+1}'*\ldots i_{s-1}'*\right)i_s'*\\
&=*i_1\ldots *i_s-i_1i_1'\ldots i_{s-1}i_{s-1}'i_s'*=i_1i_1'\ldots i_si_s'\,;\end{align*}
the last equality is just the identity to be proved in the first case, which is straightforward (again, by induction on $s$).  

If the first expansion block is $*i_1i_1'\ldots i_si_s'$ and $r_1=1$, we let 
\[\casecc{u_1}{:=}{*i_1i_1'\ldots i_{s-1}i_{s-1}'i_s'*}{s\ge 1}{-*}{s=0}\,\spa\spa v_1:=0\,.\]
Note that this definition can be obtained from Case 1 above by multiplying the corresponding terms by $q=-*$ on the left. In a similar way, we deduce the definition of $u_1$ and $v_1$ from Case 2 when $r_1=0$, and define $u_m$ and $v_m$ when the last expansion block is $i_1i_1'\ldots i_si_s'*$. Furthermore, it is again clear that, in each situation, the corresponding expansion block can be expressed as $u_1+v_1$, respectively $u_m+v_m$. 

In order to define the expansion of the whole product $\pi(A)\pi'(B)$, we need the following Lemma.
\begin{lem}\label{lem2} Consider elements $x_k=u_k+v_k$ for $k=1,\ldots,m$ in a noncommutative algebra. The following identity holds:
\[x_1\ldots x_m=v_1\ldots v_m+\sum_{\emptyset\ne I\subseteq[m]}(-1)^{|I|+1}y_1\ldots y_m\,,\]
where
\[\casec{y_k}{:=}{u_k}{k\in I}{x_k}\]
for $k=1,\ldots,m\,$.
\end{lem}
\noindent {\em Proof.} Fix a set $K$ with $\emptyset\ne K\subseteq[m]$, and let
\[\case{z_k}{:=}{u_k}{k\in K}{v_k}\]
Consider the expansion of the product $x_1\ldots x_m$ in the left-hand side of the identity to be proved. The number of times the term $z_1\ldots z_m$ appears in the expansion of the right-hand side is
\[\sum_{\emptyset\ne I\subseteq K}(-1)^{|I|+1}=1\,.\;\;\;\;\;\;\;\;\;\Box\]

Clearly, Lemma \ref{lem2} can be generalized by considering extra factors (corresponding to the intermediate blocks) between the factors in the noncommutative products above. Now let $x_k$ in this more general form of the Lemma be the expansion blocks in $\pi(A)\pi'(B)$ (with the appropriate expansions $u_k+v_k$), and expand $\pi(A)\pi'(B)$ as indicated by the Lemma.  

\underline{Step 3.2: Expanding the terms $\pi'(B)\pi(A)$.} Similar expansions are defined for the expansion blocks in such a term. Assume there are $m$ expansion blocks, and let the $k$th such block be $j_1'j_1\ldots j_t'j_t$. Based on the labels $l_k,\,r_k$ defined in Step 2 and the second relation (3) in the definition of ${\mathcal E}_{X}$, we express the considered block as $u_k+v_k$, where $u_k$ and $v_k$ are given by the rules below.
\begin{itemize}
\item Case $1'$: If $l_k=0$ and $r_k=1$ then 
\[u_k:=j_1'j_1\ldots j_{t-1}'j_{t-1}j_t*\,,\spa\spa v_k:=(-1)^t*j_1'\ldots *j_t'\,.\]
\item Case $2'$: If $l_k=1$ and $r_k=0$ then 
\[u_k:=-*j_1'j_2'j_2\ldots j_t'j_t\,,\spa\spa v_k:= j_1*\ldots j_t*\,.\]
\item Case $3'$: If $l_k=1$ and $r_k=1$ then 
\[u_k:=-*j_1'j_2'j_2\ldots j_t'j_t+j_1'j_1\ldots j_{t-1}'j_{t-1}j_t*+*j_1'j_2'j_2\ldots j_{t-1}'j_{t-1}j_t*\,,\spa\spa v_k:=0\,;\]
if $t=1$, then the third term in $u_k$ is missing.
\item Case $4'$: If $l_k=0$ and $r_k=0$ then 
\[u_k:=0\,,\spa\spa v_k:=\sum_{i=0}^t (-1)^{i} *j_1'\ldots *j_i'j_{i+1}*\ldots j_t*\,.\]
\end{itemize}
 Again, the four ways to expand a block all generalize the second relation (3). 

The whole product $\pi'(B)\pi(A)$ is expanded based on Lemma \ref{lem2}, in the same way as the product $\pi(A)\pi'(B)$. 

\underline{Step 4: Cancellation.} We claim that, when all terms in $\varSigma(X,d)$ (which we call parents) are expanded in the way shown above, the resulting terms (called sons) cancel out. Note first that there are three types of sons: 
\begin{enumerate}
\item the terms coming from the expansion of $ ...  y_1 ...  y_2  ...\;...\;...  y_m ... $ in Lemma \ref{lem2} which correspond to terms of the form $\pi(A)\pi'(B)$ in $\varSigma(X,d)$ (the intermediate blocks are indicated by $...$); 
\item the similar terms which correspond to terms of the form $\pi'(B)\pi(A)$;
\item the terms coming from the expansion of $ ...  v_1 ...  v_2  ...\;...\;...  v_m ... $ in Lemma \ref{lem2}.
\end{enumerate}
In the first class of terms we also include $*\pi(X)$ and $-\pi'(X)*$ if $d=|X|+1$. We claim the the cancellation takes place within the three classes of terms. 

\underline{Step 4.1: Cancellation in class (1).} Fix a term of the form $\pi(A)\pi'(B)$ in $\varSigma(X,d)$, and assume it has $m$ expansion blocks. Fix a set $\emptyset\ne I\subseteq[m]$ of blocks which are expanded by Lemma \ref{lem2}, meaning that the those blocks are replaced by the corresponding $u_k$. Assuming $u_k\ne 0$ for all $k$ in $I$, each of them is a sum of at most three terms $t_{k1}$, $t_{k2}$, $t_{k3}$ (cf. cases 1, 2, and 3 in Step 3). Note that, in all possible cases, the sign of a term $t_{k\cdot}$ within $u_k$ is $(-1)^{{\rm left}(t_{k\cdot})+1}$, where ${\rm left}(t_{k\cdot})$ is 1 or 0, depending on whether $t_{k\cdot}$ starts with a $*$ or not; to be more specific in the special case $u_k=\pm *$, we set ${\rm left}(*):=0$ if $*$ is at the beginning of $\pi(A)\pi'(B)$, and ${\rm left}(*):=1$ if $*$ is at the end. Upon expanding $\pi(A)\pi'(B)$, we obtain sons which are concatenations of three types of blocks: $t_{k\cdot}$ for $k\in I$, the expansion blocks which are not expanded (i.e., those indexed by $k$ in $[m]\setminus I$), and the intermediate blocks. The sign of such a son is
\begin{equation}\label{sign}(-1)^{|I|+1}(-1)^{\sum_{k\in I} ({\rm left}(t_{k\cdot})+1)}=(-1)^{1+\sum_{k\in I} {\rm left}(t_{k\cdot})}\,;\end{equation}
here we took into account the sign contribution from Lemma \ref{lem2} as well.
 
A typical term in the first class is of the form
\begin{equation}\label{word}   [...b_1']*a_1 ...  b_2'*a_2 ...\;...\;...  b_s'*[a_s...]\,,\end{equation}
where the subwords enclosed by square brackets might or might not appear. Let us count how many times this term appears and let us figure out its signs, by examining the parents $\pi(A)\pi'(B)$ in $\varSigma(X,d)$ which can give rise to it upon expansion. The $i$th symbol $*$ can combine with another letter in one of the following two ways.
\begin{itemize}
\item Case A: It combines with the letter $a_i$ after it, if any.
\item Case B: It combines with the letter $b_i'$ before it, if any.
\end{itemize}
Assuming the mentioned letters exist, the corresponding parents, written in block form, contain the adjacent pair $a_ia_i'$, respectively $b_ib_i'$. If $b_1'$ is absent, we have the following two cases.
\begin{itemize}
\item Case A: The first $*$ combines with $a_1$, as explained above.
\item Case B$'$: The first $*$ remains unchanged, as the $*$ at the beginning of a parent (see also the remarks below). 
\end{itemize}
The situation in which $a_s$ is absent is similar, leading to cases A$'$ and B. Hence, there are $2^s$ parents of the form $\pi(A)\pi'(B)$; indeed, the labeling ensures that the order of the letters is the correct one in all cases. Furthermore, according to (\ref{sign}), each of the $s$ positions in which expansion occurs has a sign contribution of $-1$ (in cases A and A$'$), respectively $1$ (in cases B and B$'$). Hence, by summing up the total sign contributions of the $2^s$ parents (obtained by taking the product of the sign contributions of the $s$ positions where expansion occurs), we get 0. 

{\em Remarks.} First, note that two consecutive $*$'s in (\ref{word}) might come from the same expansion block in the parent, which is the reason for the third term in the definition of $u_k$ in Case 3 of Step 3. Secondly, note that the terms $*\pi(X)$ and $-\pi'(X)*$, which appear with a negative sign in $\varSigma(X,d)$ when $d=|X|+1$, are also needed in the cancellation process above. Indeed, these words are special cases of (\ref{word}), cf. cases A$'$ and B$'$ above.  

\underline{Step 4.2: Cancellation in class (2).} The terms in the second class can be cancelled by a procedure similar to the one in Step 4.1. In fact, this case is easier because no expansion blocks containing $*$ are involved. 

\underline{Step 4.3: Cancellation in class (3).} We show that the terms in the third class cancel in pairs. A son of $\pi(A)\pi'(B)$ in the third class is obtained by replacing all expansion blocks $i_1i_1'\ldots i_si_s'$ with $*i_1\ldots *i_ji_{j+1}'*\ldots i_s'*$ for some $j$ with $0\le j\le s$ (cf. cases 1, 2, and 4 in Step 3). Furthermore, the labeling ensures that if $j\ge 1$, then the first $*$ is preceded by a letter $a$, and if $j<s$, then the last $*$ is followed by a letter $b'$. If there is an initial expansion block $*i_1i_1'\ldots i_si_s'$, it will be replaced by $*i_1'*\ldots i_s'*$; again, the last $*$ has to be followed by a letter $b'$. Similarly for a possible expansion block ending in $*$. Let us now replace each adjacent pair $*i'$ in the son by $i'i$ and each $j*$ by $j'j$ (cf. cases $1'$, $2'$, and $4'$ in Step 3). We obtain a word of the form $\pi'(D)\pi(C)$ with no $*$'s, up to trivial commutations given by relation (4) in the definition of ${\mathcal E}_X$; indeed, the labeling ensures that the order of the letters is the correct one. Hence $\pi(A)\pi'(B)$ and $\pi'(D)\pi(C)$ have two identical sons (which we pair up); furthermore, their signs are the same, namely $(-1)^k$, where $k$ is the number of adjacent pairs $*i'$ in the sons. Finally, it is easy to see that all sons in the third class of a typical term $\pi'(D)\pi(C)$ containing no $*$'s arise in the way mentioned above, and are paired up uniquely. \eproof

{\em Proof of the general case of Lemma {\rm \ref{lem1}}.} This follows from the special case of the Lemma. Let $X=X_1\cup X_2$, with $X_1=\{i_1,\ldots,i_s\}$, $X_2=\{j_1,\ldots,j_t\}$, and $C(i_1,\ldots,i_s,p,j_1,\ldots,j_t,q)$. Assume that both $X_1$ and $X_2$ are nonempty, because otherwise we are in the special case above. We have
\begin{equation}\label{comm}\varSigma(X,d)=\sum_{(A_1,B_1,A_2,B_2)}\pi(A_2)\pi(A_1)\pi'(B_1)\pi'(B_2)-\sum_{(C_1,D_1,C_2,D_2)}\pi'(D_1)\pi'(D_2)\pi(C_2)\pi(C_1)\,,\end{equation}
where the first summation ranges over all quadruples of sets $(A_1,B_1,A_2,B_2)$ satisfying
\begin{itemize}
\item $X_1\subseteq A_1\cup B_1\subseteq X_1\cup\{p,q\}\,,$ $A_2\cup B_2= X_2\,$;
\item $|A_1|+|B_1|+|A_2|+|B_2|=d\,$;
\item $A_1\cup A_2\ne\emptyset\ne B_1\cup B_2$ and $p\not\in A_1$, $q\not\in B_1\,$;
\end{itemize}
the second summation ranges over all quadruples of sets $(C_1,D_1,C_2,D_2)$ satisfying the conditions obtained from the ones above upon substitutions $(A_1,B_1)\rightarrow(C_2,D_2)$, $(A_2,B_2)\rightarrow(C_1,D_1)$, $X_1\leftrightarrow X_2$. 

Applying the special case of the Lemma, and assuming that $A_2\ne\emptyset\ne B_2$ are fixed, we have that
\[\sum_{(A_1,B_1)}\pi(A_1)\pi'(B_1)=\sum_{(\widetilde{A}_1,\widetilde{B}_1)}\pi'(\widetilde{B}_1)\pi(\widetilde{A}_1)\,;\]
here the first summation ranges over all $(A_1,B_1)$ satisfying the conditions above, and the second summation ranges over $(\widetilde{A}_1,\widetilde{B}_1)$ satisfying $\widetilde{A}_1\cup\widetilde{B}_1=X_1$ and $|\widetilde{A}_1|+|\widetilde{B}_1|=|A_1|+|B_1|$. Taking into account the case when $A_2=\emptyset$ or $B_2=\emptyset$, we can write the first summation in (\ref{comm}) as
\begin{align*}
&\sum_{(\widetilde{A}_1,\widetilde{B}_1,A_2,B_2)}\pi(A_2)\pi'(\widetilde{B}_1)\pi(\widetilde{A}_1)\pi'(B_2)+\casex{\pi(X_2)*\pi(X_1)-\pi'(X_1)*\pi'(X_2)}{d=|X|+1}{0}  \\
&=\sum_{(\widetilde{A}_1,\widetilde{B}_1,A_2,B_2)}\pi'(\widetilde{B}_1)\pi(A_2)\pi'(B_2)\pi(\widetilde{A}_1)+\casexsc{\pi(X_2)*\pi(X_1)-\pi'(X_1)*\pi'(X_2)}{d=|X|+1}{0}
\end{align*}
here the summations range over quadruples $(\widetilde{A}_1,\widetilde{B}_1,A_2,B_2)$ satisfying 
\begin{itemize}
\item $\widetilde{A}_1\cup\widetilde{B}_1=X_1\,,$ $A_2\cup B_2= X_2\,$;
\item $|\widetilde{A}_1|+|\widetilde{B}_1|+|A_2|+|B_2|=d\,$;
\item $\widetilde{A}_1\cup A_2\ne\emptyset\ne \widetilde{B}_1\cup B_2\,$.
\end{itemize}

By a similar procedure, the second summation in (\ref{comm}) can be written in the same way, whence $\varSigma(X,d)=0\,$. \eproof

\section{A new approach to computing the $K$-theory structure constants}

We now present the $K$-theory versions of (\ref{cor1}) and Conjecture \ref{conj1}. 

Consider the evaluation of Grothendieck polynomials at $K$-theoretic Dunkl elements
\[\fg_w(\kappa):=\fg_w(\kappa_1,\ldots,\kappa_{n-1})\,.\]

\begin{cor}\label{strconst} Given $u,\,v,\,w$ in $S_n$ with $l(w)\ge l(u)+l(v)$, we have
\[c_{uv}^w=\langle\mbox{coefficient of $w$ in $\fg_u(\kappa)\,v$}\rangle\]
in the Bruhat representation. 
\end{cor}

\begin{conj}\label{nonnk} ($K$-theoretic nonnegativity conjecture) \cite{yonpc}. For any $w$ in $S_n$, the $k$th graded piece in the evaluation
\[(-1)^{k-l(w)}\fg_w(\kappa)\]
lies in the cone ${\mathcal E}_n^+$. 
\end{conj}

\begin{rem}
{\rm 
\begin{enumerate}
\item  Corollary \ref{strconst} is a new approach to the study of the $K$-theory structure constants $c_{uv}^w$, which is superior to the one based only on Theorem \ref{kmonk}. Indeed, note first that the Littlewood-Richardson problem in $K$-theory is reduced to simplifying $\fg_u(\kappa)$, by using the relations of the quadratic algebra. Let us also note that, essentially, we are simultaneously deriving the rules for multiplying $\fg_u(x)$ by {\em all} Grothendieck polynomials. In order to better understand this fact, let us label a cover $w<wt_{ij}$ in Bruhat order by the pair $(i,j)$, and, based on Corollary \ref{strconst}, let us think of $c_{uv}^w$ as counting certain saturated chains from $v$ to $w$, which are determined by $u$. The above procedure has the advantage that we are considering the saturated chains with the same sequences of labels just once, rather than for each relevant Grothendieck polynomial $\fg_v(x)$ by which we multiply $\fg_u(x)$. The importance of Corollary \ref{strconst} is also underlied by the next remark, as well as by Conjecture \ref{nonnk}, as discussed in the fourth remark below. 
\item Working in the quadratic algebra rather than in its Bruhat representation makes sense, because its relations are simpler (the complete list of relations for the Bruhat representation is not even known), and stronger results are revealed, such as Theorem \ref{commute} and, possibly, Conjecture \ref{nonnk}. 
\item Conjecture \ref{nonnk} implies that $(-1)^{l(w)-l(u)-l(v)}c_{uv}^w\ge 0$, as proved by Brion in \cite{bripgg}. But the conjecture does not appear to follow from Brion's result. A. Yong \cite{yonpc} checked the conjecture for $S_3$.
\item If Conjecture \ref{nonnk} is true, then a combinatorial description of $\fg_u(\kappa)$ as an alternating linear combination of monomials in $[ij]$ with $i<j$ would immediately lead to a combinatorial description of the $K$-theory structure constants $c_{uv}^w$, for all $v$ and $w$. Note that a similar approach led to a proof of the Pieri-type formula for the cohomology of the flag variety in \cite{posqvp}.
\end{enumerate}
}
\end{rem}

\vspace{-1.8mm}

Fomin and Kirillov proved in \cite{fakqad} that the subalgebra of the quadratic algebra ${\mathcal E}_n$ generated by the Dunkl elements is isomorphic to $H^*(Fl_n)$. The author and A. Yong \cite{layktf} conjectured that $K^0(Fl_n)$ can be realized in a similar way, as a different commutative subalgebra of the {\em same} quadratic algebra, based on the $K$-theoretic Dunkl elements.

\begin{conj}\cite{layktf}
The ring $K(Fl_n)$ is isomorphic to the commutative
subalgebra of ${\mathcal E}_{n}$ generated by the elements $\kappa_{i}$, for 
$i=1,\ldots,n$.
\end{conj}

We have checked this conjecture for $n\leq 4$. Clearly, the main result of this paper is a step in this direction. On the other hand, we have proved other related results, including the fact that the $K$-theoretic Dunkl elements sum to 0. However, completing the proof of the conjecture will require much more effort. Indeed, its cohomology counterpart in \cite{fakqad} is already somewhat complex, while passage from cohomology to $K$-theory seems to involve a dramatic increase in complexity, as suggested by the proof of Theorem \ref{commute}.

\end{document}